\hsize=11.5cm
\vsize=19.5cm
\hoffset=2.5cm
\voffset=3cm
\nopagenumbers

\def\ni{\noindent}
\def\bn{\bigbreak\noindent}
\def\mn{\medskip\noindent}
\def\sn{\smallskip\noindent}
\def\hb{\hfil\break}
\def\ie{{\it i.e.\/}\ }
\def\eg{{\it e.g.\/}\ }
\def\PP{{\rm P}}
\def\CC{{\rm C}}
\def\ZZ{{\rm Z}}
\def\NN{{\rm N}}
\def\pn{{\rm P}^n}

\def\po{{\rm P}^1}
\def\pt{{\rm P}^3}
\def\co{{\cal O}}
\def\bqb{{\rm\bar Q}}
\def\codim{\mathop{\rm codim}}
\def\sing{\mathop{\rm sing}}
\def\chow{\mathop{\rm Chow}}
\def\pa{\partial}

\def\ve{\varepsilon}

\font\title=cmbx10 scaled\magstep2
\font\scap=cmcsc10

\vglue 4cm
\baselineskip=7mm
\ni{\title Complete intersections and rational\hb
\ni equivalence}

\baselineskip=5mm

\bn{\scap R. Barlow}

\bn{Department of Mathematical Sciences, University of Durham,\hb
\ni Durham DH1 3LE, UK}

\vskip 3cm

\baselineskip=3mm
\ni A new criterion for rational equivalence of cycles on a projective variety
over an algebraically closed field is given, and  some consequences considered.

\baselineskip=5mm
\vskip 0.8cm

\ni{\bf Introduction}

\sn The motivation for this paper was  Bloch's
conjecture [4]  on injectivity of the Abel Jacobi mappping for complex
 projective surfaces 
with $p_g=0$, a surprising conjecture in that for $p_g\neq 0$ this map
has infinite-dimensional kernel [15]. 
    This and related  conjectures (for example in [5,14]) imply the 
existence 
of many more rational equivalences under  certain  conditions 
(for example if $p_g=0$ [4], or if the variety is defined over a number field 
[5]) than otherwise. In cases where  this can be proved, the methods
rely  on special properties such as the variety having enough automorphisms.
So   a step towards finding a general method
might be to reconsider the definitions
of rational equivalence of cycles on a general variety.

The main result of this paper  (theorem {\sl2.1\/})  is a  criterion
for rational equivalence, derived from
standard definitions  ({\sl1.1\/}).
For nonsingular $V$ in $\pn$, {\sl2.1\/} states that $X$ is rationally 
equivalent to $Y$
 if and only if 
$X-Y=\bigl((A-B)\cdot D\cdot V\bigr)$, where $A$, $B$ are hypersurfaces  
of like
degree  and $D$ is a  complete-intersection cycle in a strict sense. This
generalizes the well-known fact that divisors $X$ and $Y$ are rationally
equivalent if and only if $X-Y=A-B$, where $A$ and $B$ are hypersurface
sections of like degree. It extends results of Severi and Samuel
(see section {\sl1\/}). The intersection product  used here is 
Samuel's (see section {\sl0\/}).

In  section {\sl3\/}, the description of the families of pairs of rationally 
equivalent cycles implied by {\sl2.1\/} is given. In terms of Mumford's 
description of these
sets, this amounts to replacing the Chow varieties of rational curves in
$S^nV$ (in the $0$-cycle case) in his correspondences by rational 
subvarieties which
are easy to describe. One consequence (corollary {\sl3.7\/}) is that complete 
intersection expressions  for multiples of points on $V$, like those found by
Roitman in [17] for complete intersections $V$
 with $p_g=0$ (but not nescessarily effiective), and similar to those 
sometimes implied when there are enough
 automorphisms (\eg [12]) (see {\sl 3.8.2\/}),
must always exist if all points on  $V$ are rationally equivalent to 
each other. One still needs a key for constructing these expressions for
noncomplete intersections $V$.

In section 4, there are more details for the case of surfaces over $\CC$
or $\bqb$, with some examples of conjectures equivalent to those of
[4,5,14]. Since the criterion of ({\sl2.1\/}) can be stated as existence of
solutions to a polynomial equation of sufficiently high degree
(see {\sl2.4\/}), there is some scope for 
algebraic methods or computer experiments. If one wanted to devise an algorithm
to test for rational equivalence between two points, one would need to know if 
the above degree could be bounded, which suggests connections with the problem 
of whether there are infinitely many rational points on the surface.
 
A  reason why these equations should be expected to have solutions
under the condition $p_g=0$ or when the underlying field is $\bqb$ is 
still lacking.
This is discussed in ${\sl4.3\/}$.
As for the role of $p_g$, the more geometric
versions  of the conjectures coming from ({\sl3\/}) are considered in
{\sl 4.5\/}. These require finding
bounds for the dimensions of certain ``special position'' loci, and suggest
 links with  vector bundles, K-theory [4,14,16], or instantons [8].
 The  calculations  seem  difficult.
 
In {\sl4.6\/} we give an example of describing rational equivalences 
explicitly. The families of
 pairs of points $X,Y$ on a generic surface in $ \pt$ which are made 
rationally equivalent  by $X-Y$ being a difference of two intersections 
 with lines are 
described. It might be interesting to find the corresponding families 
obtained by intersecting with pairs of complete intersections of higher 
degree. 
\headline{\hfil{\scap R.  Barlow}\hskip4.3cm\tenrm\folio}
\pageno=2

   Not pursued here is the possibility of using the construction
of section {\sl3\/}  backwards, \ie
to obtain information about Chow varieties of varieties such as $\pn$
known to have trivial Chow groups.

\bigbreak
\ni{\bf0. Preliminaries on cycles and intersections.}

\sn For a projective variety $V\subset\pn$ over an algebraically closed field
$k$, the $r$-cycles (formal sums of irreducible subvarieties of dimension $r$)
form a group $Z_r(V)$. The $r$-cycles rationally equivalent (written $\sim$)
to zero form a subgroup $B_r(V)$; and the quotient
$A_r(V)=Z_r(V)/B_r(V)$, also called $CH_r(V)$, is the Chow group.

The degree of a cycle $\sum n_i X_i$ on $V$ is $\sum n_i \deg X_i$.
The classes of degree zero cycles form a subgroup $A^0_r(V)$. The distinct
cycles $X=\sum n_i X_i$ of degree $l$, with $n_i\geq0$ for all $i$, are
in one-to-one correspondence with the points of a projective algebraic set
${\chow}^l_r(V)$ (see [18]). In particular for $r=0$,
${\chow}^l_0(V)$ is isomorphic to the symmetric product $S^lV$.

If $W$ is a subscheme of $V$, then $[W]$ denotes the associated cycle
$\sum l(\co_{W_i,W})W_i$ (see [9]), where the $W_i$ are the irreducible
components of $W$ and $l(\co_{W_i,W})$ is the length of the Artinian ring
$\co_{W_i,W}$. When $W$ is a subvariety, we write $[W]=W$.

If $X$ and $Y$ are cycles meeting properly on $V$, \ie such that
$\codim X\cap Y = \codim X + \codim Y$, Samuel's intersection cycle
$(X\cdot Y)_V$ is a well-defined member of $Z(V)$ when $V$ is nonsingular
at generic points of $X\cap Y$.
It agrees as a formal sum with Fulton's more generally defined product
$\overline X \cdot\overline Y$ in $A_*(X\cap Y)$ ([9] chapter 7), and suffices 
for this
paper. In fact, we only use it for the case when at least one of $X$ or $Y$
is the cycle associated to a local complete intersection scheme in $V$, when
the multiplicities of components are again lengths of local rings. So if
$X$ and $Y$ are cycles associated to schemes $X_0$ and $Y_0$ which meet
properly, and $Y_0$ is a local complete intersection in $V$ with $V$
nonsingular along $X_0\cap Y_0$, we have $(X\cdot Y)_V=[X_0\cap Y_0]$.
When $V={\pn}$ we usually drop the suffix $V$.

For example, if $f$ is a rational function on $\pn$ represented by $a/b$,
then the divisor $(f)$ of $f$ is defined by
$(f) = \bigl[\{a=0\}\bigr] - \bigl[\{b=0\}\bigr]$, and the divisor $(f_W)$ 
of $f$ restricted to
a variety $W\subset \pn$ (on which it is defined and nonzero) could be defined
by $(f_W) = \bigl((f)\cdot W\bigr)_{\pn}$. This was Samuel's definition of 
$(f_W)$, and it agrees with the $[{\rm div}\, f_W]$ in chapter 1 of [9].
If $W\subset V\subset\pn$, we also have $(f_W)=\bigl((f_V)\cdot W\bigr)_V$ 
if $V$ is
nonsingular at generic points of $(f_W)$, but not always (see [18], page 103).

\eject

\ni{\bf1. Some characterizations of rational equivalence.}

\sn{\sl1.1  Definitions}
\sn Let $X$ and $Y$ be $r$-cycles on $V$. Then $X$
is rationally equivalent to $Y$ on $V$, written $X\sim Y$, if and only if one
of the following three conditions holds. The equivalence of these three,
at least for nonsingular $V$, is proved in [19].
Samuel used (2) as a definition; (1) and (3) are used in [9] and [15]
respectively, for example.

\sn
\item{(1)} $X-Y = \sum_{i=1}^t (f_{D_i})$, where $D_1,\ldots,D_t$ are
  subvarieties of dimension $r+1$ and $f_{D_i}$ is a nonzero rational function
on $D_i$. Note that the $D_i$ can be assumed distinct.
\item{(2)} $X-Y = \pi_*\bigl(V\times\{0\}\cdot D\bigr)- 
  \pi_*\bigl(V\times\{\infty\}\cdot D\bigr)$,
where $D$ is an $(r+1)$-cycle on $V\times\po$ and $\pi$ is projection onto $V$.
\item{(3)} (Assuming that $X$ and $Y$ are effective of degree $l$ and
  identifying them with their Chow forms) there exists an integer $e$ and a
  morphism $f:\po\to\chow^{l+e}_r(V)$ such that $f(0)=X+Z$ and
  $f(\infty)=Y+Z$ for some $Z$ in ${\chow}^e_r(V)$.

\sn
Earlier, \lq\lq rational equivalence\rq\rq\ between $X$ and $Y$ was defined by
Severi [21, 2] as \lq\lq belonging to the same intersection series\rq\rq\
(\ie $X=\sum m_i X_i$ and $Y=\sum m_i X_i'$, where
$X_i=H_{i,1}\cdot\ldots\cdot H_{i,m-r}$ and
$X_i'=H_{i,1}'\cdot\ldots\cdot H_{i,m-r}'$ with divisors
$H_{i,j}\sim H_{i,j}'$ for all $i,j$). This was shown to be equivalent to
insisting that $X-Y=\sum m_i(X_i-X_i')$ with $X_i$ and $X_i'$ as above, by
[23]. (See [9] for an historical outline.)

These definitions were shown to be equivalent to a version of (2) above by
Baldassari [2], using his generalization of {\sl1.2\/}\ (1) below to 
cycles on $V\times\pn$, and {\sl1.2\/}\ (2).

\mn{\sl1.2 Severi's theorems}
\sn(1) (See [18].) {\sl Every $r$-cycle $X$ on $\pn$ can be written as
$X=H_1\cdot\ldots\cdot H_{n-r}$, where the $H_i$ are  divisors
(in general not effective).}\hb
\ni(2)(see [21]). {\sl On $V\subset\pn$ nonsingular, $X\sim Y$ if and only if
 $X$ and $Y$
\lq\lq belong to the same series of intersection away from some fixed and
semi-fixed components\rq\rq.}

\ni Some related work of Severi is discussed  in [27]. 
  It was assumed that $A_0^0$ is always
finite-dimensional; this was disproved by Mumford [15]. 

\mn{\sl1.3 Samuel's theorem}
\sn In theorem 10 of [19], Samuel proved the following
natural generalization of linear equivalence of divisors.

\ni(1) Theorem. {\sl On a nonsingular variety $V\subset\pn$, $r$-cycles $X$ and
$Y$ are rationally equivalent if and only if $X-Y=\bigl((f)\cdot D\bigr)$ 
for some rational function $f$ on $\pn$ and $r$-cycle $D$ on $V$.}

\ni To prove this, after showing that his definition ((2) above) implied (1),
Samuel used the following trick to replace the sum of terms
$\sum_{i=1}^t (f_{D_i})$ or $\sum_{i=1}^t\bigl((f_i)\cdot D_i\bigr)$ (see the 
previous section) by a single expression, with $D=\sum D_i$.

\ni(2) Trick. {\sl Let $f_i$ be a rational function on $\pn$ restricting to
$f_{D_i}$ on $D_i$ for $i=1,\ldots,t$. Let $g_i$ be a form vanishing on every
$D_j$ $(j\neq i)$ but not on $D_i$. Choose $g_1,\ldots,g_t$ all of the same
(high enough) degree. Let $f=\sum g_i f_i/\sum g_i$. Then $f=f_i$ on $D_i$.}

\ni This requires $D_i\neq D_j$ for $i\neq j$, which can be assumed in
definition (1) by
combining: $(f_{D_i}) + (f_{D_i}') = (f_{D_i}f_{D_i}')$.

\bigbreak
\ni{\bf2. Complete intersection expression for rational equivalence.}

\sn We prove {\sl2.1\/}, which is like a version of Severi's theorem
{\sl1.2} (2) with hypersurface sections (in particular effective divisors)
replacing arbitrary virtual divisors. Severi also stated in [20] that
effective divisors would suffice for 0-cycles on surfaces (but deferred
publishing the complete proof). The proof of {\sl2.1\/} given here amounts 
to showing that
in Samuel's theorem {\sl1.3\/} we can choose $D$ to be a complete
intersection cycle on $V$. The condition of {\sl2.1\/} can be expressed
compactly as an equation on Chow forms (Corollary {\sl2.4\/}). There is some
choice in the expression (see {\sl2.5\/}), and we may try to simplify it by
re-embedding the variety (see {\sl2.6\/}).

\mn{\sl2.1 Theorem}
\sn{\sl Let $V$ be a variety of dimension $m$ in $\pn$, nonsingular away from
a set of dimension $r$, where $r\leq m-2$. Then $r$-cycles $X$ and $Y$ on $V$
are rationally equivalent on $V$ if and only if $X-Y = V_X-V_Y$, where
$V_X=\bigl[\{a=h_1=\ldots=h_{m-r-1}=0\}\cap V\bigr]$ and
$V_Y=\bigl[\{b=h_1=\ldots=h_{m-r-1}=0\}\cap V\bigr]$
for some forms $a$ and $b$ of like degree and forms $h_1,\ldots,h_{m-r-1}$
such that $V_X$ and $V_Y$ are in $Z_r(V)$.}

\sn Remarks.
(1) We could also characterize rational equivalence of $X$ and $Y$
modulo some fixed set $Z \subset V$ in this way, \ie for $X$ and $Y$
in $Z_r(V\backslash Z)$, $X-Y=V_X-V_Y +\ve$, with $\ve$ supported on $Z$,
and for this could allow $V$ to be singular on $Z$ of higher dimension.\hb
\ni(2) If $V$ is nonsingular in dimension $r$, we can write $V_X=(A\cdot D)_V$
and $V_Y=(B\cdot D)_V$, where $A=\bigl[\{a=0\}\cap V\bigr]$, 
$B=\bigl[\{b=0\}\cap V\bigr]$ 
and $D=\bigl[\{h_1=\ldots=h_{m-r-1}=0\}\cap V\bigr]$, as in {\sl0\sl}.

\sn Proof. The sufficiency of the conditions is clear. So assume $X\sim Y$.
By definition (1), $X-Y=\sum_{i=1}^{t_0} (f_{D_i})$ where $D_1,\ldots,D_{t_0}$
are distinct $(r-1)$-dimensional subvarieties.
Assuming Proposition {\sl2.2\/} below, we can find a complete intersection
cycle $D=\bigl[\{h_1=\ldots=h_{m-r-1}=0\}\cap V\bigr]$ such that
$D=\sum_{i=1}^t D_i$ for $t\geq t_0$, and $D_i\neq D_j$ for $i\neq j$.
Let $\lambda_{t_0+1}, \ldots,\lambda_t$ be arbitrary nonzero constants
and set $f_{D_i} = \lambda_i$ for $i\geq t_0+1$.
Then since $(f_{D_i})=0$, we still have $X-Y=\sum_{i=1}^t (f_{D_i})$.
Now apply Samuel's trick {\sl1.3\/}\~(2). Let $f=\sum g_i f_i /\sum g_i$,
where $f_i$ is an extension of $f_{D_i}$ to $\pn$ and $g_i$ is a form
vanishing on every $D_j$ ($j\neq i$) but not on $D_i$.
Then $(f_{D_i})=\bigl((f)\cdot D_i\bigr)_{\pn}$, 
so $X-Y=\bigl((f)\cdot D\bigr)_{\pn}$.
Let $a$ and $b$ be forms such that $a/b$ represents $f$ in lowest terms.
Then $V_X=\bigl[\{a=h_1=\ldots=h_{m-r-1}=0\}\cap V\bigr]$ and
$V_Y=\bigl[\{b=h_1=\ldots=h_{m-r-1}=0\}\cap V\bigr]$ gives {\sl2.1\/}.

\mn{\sl2.2 Proposition}
\sn{\sl Let $V\subset\pn$ be an algebraic set of dimension $m$; 
and let $D_0$ be an $l$-cycle on $V$ with $l\leq m-1$, such that
$D_0=\sum_{i=1}^{t_0} D_i$ with $D_i$ irreducible and $D_i\neq D_j$
for $i\neq j$. Suppose that $V$ is nonsingular away from at most a set $Z$
of dimension $l-1$. Then there exist forms $h_1,\ldots,h_{m-l}$ on $\pn$
such that $\bigl[\{h_1=\ldots=h_{m-l}=0\}\cap V\bigr]=D$ 
where $D=\sum_{i=1}^t D_i$ for $t\geq t_0$ and $D_i\neq D_j$ ($i\neq j$).}\hb
\ni To prove this, we use the following Lemma.

\mn{\sl2.3 Lemma}
\sn{\sl Let $V$, $D_0$ and $Z$ be as in 2.2. Then there exists a form $h$
vanishing on $D_0$ such that $\{h=0\}\cap V=V'$, where $V'$ has dimension
$m-1$ and is nonsingular away from a set $Z'$, still of dimension at most
$l-1$.}
\sn Proof. Suppose $V=\cup_{i=1}^s V_i$ where $V_i$ is irreducible.
We show below that there exist forms $g_i$ (for $i=1\ldots s$),
 all of some high enough degree $d$, vanishing
on $\bigl(D_0\cup(\cup_{j\neq i} V_j)\bigr)=X_i$, such that 
$\{g_i=0\}\cap V_i$ is
nonsingular on $V_i$  not only away from $Z$ and $D_0$, but also 
 at generic points of $D_0$. Then we can take
$h=\sum g_i$ to give {\sl2.3\/}.

By Hilbert's basis theorem, the base locus of the linear system $I_d(X_i)$,
(forms of degree $d$ vanishing on $X_i$), is just $X_i$ for large
enough $d$. By Bertini's theorem II ([1]) , generic
$g_i\in I_d(X_i)$ has $\{g_i=0\}\cap V_i$ nonsingular away from
$X_i\cup\sing V_i$, \ie  away from   $D_0\cup\sing V_i$.

At a generic point $p\in D_0\cap V_i$, both $D_0$ and $V_i$ are
nonsingular. Since $\dim D_0<\dim V_i$, we have a strict inclusion
$\bigl\{\pa f\big|_p:f\in I_d(V_i)\bigr\}
\subset\bigl\{\pa f\big|_p :f\in I_d(X_i)\bigr\}$
(when $d$ is large enough). So for generic $g_i\in I_d(X_i)$,
$\{g_i=0\}\cap V_i$ is nonsingular at $p$.

\ni To prove {\sl2.2\/}, we apply {\sl2.3\/} $m-l$ times in succession.

\mn{\sl2.4 Corollary of 2.1}
\sn{\sl Let $\phi_X$ denote the Chow form for a cycle $X$ (see [18]). Then in
the notation of {\sl2.1\/}, $X\sim Y$ on $V$ if and only if
$\phi_X \phi_{G_1} = \phi_Y \phi_{G_2}$ for some r-cycles
$G_1$ and $G_2$ both  cut out on $V$ by complete intersections of the  same
multidegrees.}

\sn The nontrivial solutions for the resulting system of equations in
$\bigl({\chow}^l_r(V)$
\break
$\times\Pi^{m-r-1}_{i=0}{\PP^{n_i}}\bigr)^2$ 
would give information about
$A^0_r(V)$ (or about $\chow^l_r(V)$, if $A^0_r(V)=0$ as for $V=\pn$).
For $r\neq 0$, this may be more accessible than  it is in {\sl2.1\/}.
The $\phi_{G_i}$ are simpler if $V$ is a complete intersection.
That is one reason for allowing $V$ to be singular above and considering
new embeddings ({\sl2.5--2.6\/}). We do not address the question of how 
many expressions $X-Y=\sum(f_{D_i})$ with a given value of 
$\sum\deg f_{D_i}$ there may be, but rather ask how to build new complete 
intersection expressions for $X-Y$ out of given ones.

\mn{\sl2.5 Varying expressions}
\sn If $V$ is nonsingular in dimension $r$ and $X-Y=(A\cdot D)_V-(B\cdot D)_V$
with $A=\bigl[\{a=0\}\cap V\bigr]$, $B=\bigl[\{b=0\}\cap V\bigr]$,
$D=\bigl[\{h_1=\ldots=h_{m-r-1}=0\}\cap V\bigr]$ as in {\sl2.1\/}, 
we can choose new
$A$, $B$, $D$ as follows. The corresponding changes in $a$, $b$, $h_i$
give similar new expressions $V_X-V_Y$ in the general case.
\sn(1) We can choose $A$ and $B$ with $A\cap B$ proper. This is proved
below.\hb
\ni(2) $X-Y=\bigl((A+E)\cdot D\bigr)_V-\bigl((B+E)\cdot D\bigr)_V$ for any 
divisor $E$ meeting $D$ properly.\hb
\ni(3) If $A\cap B$ is proper and $\bigl[\{h_i=0\}\cap V\bigr]=H_i$, so that
$D=(H_1\cdot\ldots\cdot H_{m-r-1})_V$, then
$X-Y=\bigl((A-B)\cdot(H_1\cdot\ldots\cdot H_j\cdot(H_j+rC)\cdot H_{j+1}\cdot
\ldots\cdot H_{m-r-1})\bigr)_V$
for generic $C$ in the pencil $\langle A,B\rangle$, for any $r\in{\NN}$.\hb
\ni(4) Applying (3) to each $H_i$ and then applying (2), we can construct new
$A$, $B$, $H_i$ for $i=1,\ldots,m-r-1$  all of the
same sufficiently high degree.
\sn To prove (1), suppose that $X-Y=\bigl((f_0)\cdot D\bigr)_V$, where
$f_0$ is given in lowest terms by $a_0/b_0$ with $\{a_0=b_0=0\}\cap D$
improper. Let $a=\gamma a_0$, $b=\gamma b_0+\delta$, where $\gamma$ is any
form of degree $d-d_0$ not vanishing on any component of $\{a_0=0\}\cap V$ or
$D$. Let $\delta$ be in $I_d(D)$. Then for sufficiently large $d$,
$\{a=b=0\}\cap V$ is proper for generic $\delta$, and $f=a/b$ restricted
to $D$ is the same as $f_0$.

\mn{\sl2.6 Multiple embeddings and projections}
\sn(1) Linear rational equivalences. By (4) above, whenever $X\sim Y$ on $V$,
there is an embedding $V\subset\pn$ (some $s$-tuple of our original embedding)
such that $X-Y=(L_1\cdot V)_{\pn} -(L_2\cdot V)_{\pn}$ with $L_1$, $L_2$ linear
subspaces of $\pn$ of codimension $m-r$, and $\dim L_1\cap L_2 =\dim L_1-1$.
The existence of such an expression for $X-Y$ on general $V\subset\pn$
(so with $s=1$ above) would put some conditions on $V$ ({\it cf. 4.6\/}),
but it is not obvious how these conditions meet on the set of $s$-tuply
embedded $V$.\hb
\ni(2) Projections. Another way to try to simplify a search for rational
equivalences would be to project $V\subset\pn$ into a smaller space
${\PP^{n'}}$ ($n>n'$). Not all rational equivalences on the image $V'$
lift to $V$ (if $V\to V'$ is not an isomorphism), as the example below 
illustrates. On the other hand, if $V\to V'$ is a birational morphism of 
surfaces (for example) then $A_0^0(V)$ is finite-dimensional (see [4])
if and only if $A_0^0(V')$ is.

\mn{\sl2.7 Example}
\sn Let $X$ and $Y$ be any points on a surface $V$ in $\pn$. Then there is a 
birational map $\phi: V\to V'\subset{\PP^3}$ such that
$\phi(X)-\phi(Y)=\bigl((L_1-L_2)\cdot V'\bigr)_{\PP^3}$ where $L_1$ and
$L_2$ are lines.

\ni Proof. Choose hypersurface sections $A$, $B$, $D$ with $(A\cdot D)=X+X_1$,
$(B\cdot D)=Y+Y_1$, for some $X_1$, $Y_1$ disjoint from $X$, $Y$. Then choose
$E$ with $(E\cdot D)=X_1+Y_1+W$ disjoint from $X+Y$. Adjusting the
$\bigl(((A+E)-(B+E))\cdot D\bigr)$ as in {\sl2.5\/} we obtain $A'$, $B'$,
$D'$ all of the same degree $s$ with $(A'\cdot D')=X+2X_1+Y_1+W'$ and
$(B'\cdot D')=Y+X_1+2Y_1+W'$. Then let $\phi$ be $s$-tuple embedding $\psi$
followed by projection to ${\PP^3}$ from a centre in $\psi(\{a'=b'=d'=0\})$,
where $a'$ defines $A'$ etc.

\bigbreak
\ni{\bf3. Families of pairs of rationally equivalent cycles.}

\sn Let $V$ be a projective variety of dimension $m$, with a chosen embedding
$V\subset\pn$ of degree $d_0$. Let $U^l$ be the subset of
$\bigl({\chow}^l_r(V)\bigr)^2$ consisting of pairs of effective cycles 
which are rationally equivalent.
          
If $V$ has singular locus of dimension at most $r$, theorem {\sl2.1\/} and
the following constructions hold. If $V$ has a singular locus $S$ of
dimension greater than $r$, we could consider pairs $X$, $Y$ of cycles
on $V\backslash S$ --- still allowing rational equivalences with auxiliary
cycles $Z\subset S$ --- and obtain similar results.


\mn{\sl3.1 Definitions}
\sn Let ${\PP^{n_i}}$ be the projective space of forms $h_i$ on $\pn$
of degree $e_i$, for $i=0,\ldots,m-r-1$.
Call $e_0=s$ and $h_0=a$.
Let $e=(e_1,\ldots,e_{m-r-1})$ be the multidegree of the complete
intersection $\{h_1=\ldots=h_{m-r-1}=0\}$, and let
$d=d_0 s \Pi_{i=1}^{m-r-1} e_i$ be the degree of
$\bigl[\{a=h_1=\ldots=h_{m-r-1}=0\}\cap V\bigr]$.
Let $W^l_{se}$ be the incidence correspondence in
${\chow}^l_r(V) \times{\chow}^{d-l}_r(V) 
\times \Pi^{m-r-1}_{i=0}{\PP^{n_i}}$, consisting of the closure of the set
of all $(X,Z,a,h_1,\ldots,
\break
h_{m-r-1})$ such that
$\bigl[\{a=h_1=\ldots=h_{m-r-1}=0\}\cap V\bigr]=X+Z$.
Let $\Delta^{d-l}_{se} \subset W^l_{se}\times W^l_{se}$ be the set of pairs 
with matching ${\chow}^{d-l}_r(V)$ and $h_1,\ldots,h_{m-r-1}$ terms. This is
a closed set.

\mn{\sl3.2 Theorem
\sn With the above notation, 
let $\pi:(W^l_{se})^2\to\bigl({\chow}^l_r(V)\bigr)^2$ be projection. 
Let $U^l_{se}=\pi(\Delta^{d-l}_{se})$ (a closed set).
Then $U^l=\cup_{s,e} U^l_{se}$.}

\sn Remarks. (1) $U^l$ is also expressed as a countable union of closed sets,
in Lemma 3 of [15]: we  give an alternative description of these sets. 
\ni(2) The sets $U^l_{se}$ belong to $r$-cycles on $V$ with its chosen
embedding $\phi:V\to\pn$, so if there is any ambiguity, one could 
label them $U^l_{se}(r,V,\phi)$. If $\{U^l_i\}_{i\in I}$ is the
collection of all irreducible components of all of the $U^l_{se}$, then
$U^l=\cup_I U^l_i$ is independent of $\phi$.\hb
\ni(3) The set $U^l_{se}$ includes the diagonal $\Delta$ of
$\bigl({\chow}^l_r(V)\bigr)^2$ (for large enough $s$, $e$). Also there may
be \lq\lq improper\rq\rq\ components on which the rational equivalence is not
given as a difference of proper complete intersections as in {\sl2.1\/}, but 
only
as a limit of such. If we let $V^l_{se}$ be the set of components of $U^l_{se}$
for which some (and hence generic, by remark below) $(X,Y)$ in $V^l_{se}$
has a proper
expression as in {\sl2.1\/}, we have
$U^l = \cup_{s,e} V^l_{se}
= \cup_{s,e} \overline{\bigl(V^l_{se}\backslash\Delta\bigr)}\cup\Delta$.

\sn Proof of {\sl3.2\/}. By {\sl2.1\/}, we know that 
$U^l\subseteq \cup_{s,e} U^l_{se}$. To show that $U^l_{se}\subset U$, we
consider the correspondence for all complete intersections of type $s,e$,
\ie $W^d_{se}$. Since
$W^d_{se}$ maps birationally onto $\Pi_{i=0}^{m-r-1} {\PP^{n_i}}$,
its image $J$ in ${\chow}^d (V)$ is at least unirational.
So any two points on $J$ are joined by a rational curve  (Lemma 4 of [19]).
If $(X,Y)\in U^l_{se}$, there is some $Z$ such that $X+Z$ and $Y+Z$ are in
$J$. By definition {\sl1.2\/}\ (3), $X\sim Y$.

\ni Remark. If we label the improper part of $W^l_{se}$ as $T$, in other 
words
$T$ is the set of all $(X,Z,a,h_1,\ldots,h_{m-r-1})$ in $W^l_{se}$ such that
$\{a=h_1=\ldots=h_{m-r-1}=0\}\cap V$ has a component of dimension $r+1$,
then $T$ is closed. We could let $\Delta_0$ be the union of the components of
$\Delta^{d-l}_{se}$ not contained in $(T\times W)\cup(W\times T)$ and let
$V^l_{se}=\pi(\Delta_0)$  is as in {\sl 3.1\/}\ (3).

Often it is convenient to know that we only need to look at one of the sets
$U^l_{se}$. This follows from the next proposition.

\mn{\sl3.3 Proposition
\sn Let $P=\cup^q_{\alpha=1} U^l_{s(\alpha),e(\alpha)}$, where
$\bigl(s(\alpha),e(\alpha)\bigr)=
  \bigl(s(\alpha),e_1(\alpha),\ldots,e_{m-r-1}(\alpha)\bigr)\break
\in{\NN^{m-r}}$.
Then there exists $(s,e)\in{\NN^{m-r}}$ such that $P\subset U^l_{se}$.}

\sn Proof. It is clear from {\sl2.4\/} that the $U^l_{se}$ will form a system 
of nests. We may assume that $q=2$. If $U^l_{s(\alpha),e(\alpha)}\subset
U^l_{s'(\alpha),e'(\alpha)}$ for $\alpha=1,2$ with $s(\alpha)=s'(\alpha)$
or $e_i(\alpha)=e_i'(\alpha)$, we say the $s$ (respectively $e_i$) term
can be matched. Rules (2) and (3) from {\sl2.5\/} become
$$\eqalign{
     (2)\qquad U^l_{s,e} &\subset U^l_{s+t,e}\quad({\rm any}\ t\in{\NN})\cr
     (3)\qquad U^l_{s,e} &\subset U^l_{s,e'}\cr}
$$
where $e_i'=e_i$ for $i\neq j$ and $e_j'=e_j+rs$ for some $r\in{\NN}$.
(In different notation, {\sl2.5\/} still holds for $V$ singular away from
$\dim r$.) It is easy to match the $s$ terms using (2), so we do this last.
To match the $e_j$ term, one finds (3) alone may not suffice, but (2)
followed by (3) does. This requires two pairs of positive integers
$t(\alpha)$ and $r(\alpha)$, $\alpha=1,2$, satisfying
$$\eqalign{
  e_j(1)+r(1)\bigl(s(1)+t(1)\bigr) &= e_j(2)+r(2)\bigl(s(2)+t(2)\bigr) \cr
                                   &= e_j'. \cr
}$$
There are many solutions. For example, assuming $e_j(1)>e_j(2)$, choose
$t(\alpha)$ ($\alpha=1,2$) such that $s(\alpha)+t(\alpha)=p(\alpha)$,
where $p(\alpha)$ ($\alpha=1,2$) are distinct primes for which the
expression $\lambda p(1)+\mu p(2)=1$ has integer solutions $\lambda>0$
and $\mu<0$. Then let $r(1)=\lambda\bigl(e_j(2)-e_j(1)\bigr)$ and
$r(2)=-\mu\bigl(e_j(2)-e_j(1)\bigr)$.

\mn{\sl3.4 Corollary
\sn If $P$ is a closed subset of $U^l$ and $k$ is uncountable, then
$P\subset U^l_{se}$ for some $(s,e)\in{\NN^{m-r}}$.}

\mn{\sl3.5 Corollary
\sn A variety $V\subset\pn$ has has $A^0_r(V)=0$ if and only if
$\bigl({\chow}^l_r(V)\bigr)^2=\cup_{s,e} U^l_{se}$. Over uncountable
fields $k$, this holds if and only if
$\bigl({\chow}^l_r(V)\bigr)^2= U^l_{se}$ for some $s$, $e$.}\hb
\ni This gives, for example, the following.

\mn{\sl3.6 Corollary
\sn A nonsingular variety $V\subset\pn$ has $A^0_0(V)=0$ if and only if
given any  zero-cycle $Z$ of degree $N$ on $V$,and any point $X$ on $V$,
there exist complete intersection cycles $D_i$ and hypersurface sections
 $A_i$, $B_i$
(of equal degree) giving $NX=\sum_{i=1}^t\bigl((A_i-B_i)\cdot D_i\bigr)+Z$.

Here $t\geq N$, with equality if $Z$ is effective. Furthermore over
uncountable $k$, we may assume that $A_i$ and $D_i$ have the same degrees
for all $i$.}

\sn Proof. By {\sl3.5\/}, $A^0_0(V)=0$ if and only if
$V\times V=U^1_{se}$ for some $s$, $e$. In fact, using {\sl3.4\/} we can
also assume $s$, $e$ large enough that {\sl proper\/} expressions
as in {\sl2.1\/} for $X\sim Y$ are given by $U^l_{se}$ for all $(X,Y)$.
For the converse, we use the fact that the kernel of the Albanese mapping
is divisible [14].

\eject

\ni{\sl3.7 Corollary
\sn A variety $V\subset\pn$ over {\CC} has $A^0_0(V)=0$ if and only if 
there exists $N\in{\NN}$ such that for each point $X$ on $V$ there is a 
cycle $E_X$ in $\pn$ with $(E_X\cdot V)=NX$.}\hb
\ni(By Severi's theorem {\sl1.2\/}\ (1), $E_X$ is an intersection of virtual
divisors, so a difference of two sums of intersections of 
hypersurfaces.)

\ni Proof. Let $Z$ be a complete intersection cycle in {\sl3.6\/}, to get
necessity. If for any point $X$ , $NX=(E_X\cdot V)$ for some $E_X$ , then
$N\delta\sim0$ for all $\delta$ in $Z^0_0(V)$, so $A^0_0(V)=0$ .

\ni Question. How is  {\sl3.7\/} related  to the result from [16]:  for an 
affine surface $V$, $A^0_0(V)=0$ if and only if every point is a complete 
intersection?

\mn{\sl3.8 Further questions}
\sn (1) Roitman [17] proved that a complete intersection $V$ with $p_g=0$
has $A^0_0(V)=0$, by constructing a cone $E_X$ with $(E_X\cdot V)=NX$
for each $X$. This $E_X$ is defined by the parts $f_{ij}$ of the Taylor
expansions at $X$, $f_i=\sum f_{ij}$, of all the equations $f_i$ for $V$
(see [14]) with some extra equations to correct the dimension 
if necessary.
 Are complete intersections $V$ are the only varieties for which
$E_X$ can be assumed to be effective in {\sl3.7\/}?
How can the $E_X$ in {\sl3.7\/} be found  in general?

\sn(2) Expressions $NX=(E_X\cdot V)+Z$,   with
$Z$ in some fixed locus, can be found for a surface $V$ with $p_g=q=0$ and
\lq\lq enough automorphisms\rq\rq\ (see [3]), if in addition the
relevant quotients $V/H_i$ by subgroups $H_i$ of
${\mathop{\rm Aut}} V$ are rational.
 The original example with enough automorphisms was the Godeaux surface 
$V=Q/{\ZZ_5}$, where $Q$ is the Fermat quintic
(see [12,4]). Here  $NX$ is expressed as a combination of fibres of maps
$V\to V_i=V/H_i$ with $H_i<{\mathop{\rm Aut}} V$. Such fibres can be written
as complete intersections minus points on a fixed divisor, when the surfaces
$V_i$ are rational, as in this case.
Can these expressions for this special  $V$
be used to  give  any information on the solutions
for generic Godeaux surface $V$, now known to have $A_0^0=0$ by other
methods [25] ?


\mn{\sl3.9 Families of varieties}
\sn Let $V\to T$ be a morphism of projective varieties such that for generic
$t\in T$, $V_t$ is a nonsingular variety of dimension $m$. Suppose that
$V\subset\pn$,so that  $V_t$ is  a subvariety of $\pn$ of degree $d$, for
generic $t$.
Let $W_{se}(T)$ be the collection of the $W_{se}$ for each $V_t$,
\ie  $W_{se}(T)$ is the closure of the set of all
$(X,Z,a,h_1,\ldots,h_{m-r-1},t)$ such that
$\bigl[\{a=h_1=\ldots=h_{m-r-1}=0\}\cap V_t\bigr]=X+Z$, and is a subset of
${\chow}^l_r(V)\times{\chow}^{d-l}_r(V)\times\Pi_{i=0}^{m-r-1}{\PP^{n_i}}
  \times T$.
Then let $\Delta^{d-l}_{se}(T)$ be the set of pairs in 
$\bigl(W_{se}(T)\bigr)^2$ with $Z$, $h_i$ (all $i$) and $t$ terms matched.
Let $\pi:\Delta^{d-l}_{se}(T)\to{\chow}^l_r(V)^2\times T$ be
projection, and let $U^l_{se}(T)=\pi(\Delta^{d-l}_{se})$.
Then if $U^l(T)$ is the set of all $(X,Y,t)$ such that
$(X,Y)\in{\chow}^l_r(V_t)^2$ and $X\sim Y$ on $V_t$, we have
$U^l(T)=\cup_{s,e} U^l_{se}(T)$.

\mn{\sl3.10 Examples of families}
\sn(1) $r$-cycles as families of $0$-cycles: 

If $X$ and $Y$ are $r$-cycles on $\pn$ of degree $l$ (or more
generally on $V\subset\pn$), we could view $X$ as a family of 0-cycles
$(X\cdot\Lambda)$ where $\Lambda\in{\rm Gr}(n-r+1,n+1)$, on the sections 
$\Lambda$ of $\pn$, and describe the triples $(A,B,D)$ ($A$, $B$ hypersurfaces
and $D$ a complete intersection in $\pn$) giving rise to rational equivalences
$X-Y=(A\cdot D)-(B\cdot D)$ in terms of correspondences $\Delta^{d-l}_{se}(T)$
for 0-cycles on $\pn$ (which are easier to describe explicitly) using the
following.

\ni{\sl Lemma. $X-Y=(A\cdot D)-(B\cdot D)$ if and only if
$(X\cdot\Lambda)-(Y\cdot\Lambda)=\bigl(((A-B)\cdot D)\cdot \Lambda\bigr)$
for generic $\Lambda$ meeting $X$ and $Y$ properly.}

\ni(2) Families of surfaces:

If $V\to T$ is a family of surfaces $V_t\subset\pn$ with
$A^0_0(V_t)=0$ for generic $t$, then there exists $(s,e)\in{\NN^2}$
such that given points $X$, $Y$ on $V_t$ there exist hypersurface sections
$A$, $B$, $D$ of $V_t$ defined by equations of degree $s$, $s$, $e$
respectively, such that $X-Y=(A\cdot D)_{V_t}-(B\cdot D)_{V_t}$.
In other words, $V_t\times V_t=V_{se}^1$ (for  generic $V_t\subset\pn$)
in the notation of {\sl3.2}\ (3).

\bigbreak

\ni{\bf4. 0-cycles on surfaces defined over $\CC$ or a number field.}

\sn{\sl4.1 Background: conjectures A and B}
\sn For more details see [4, 14].
   Let F be a complex projective surface. Mumford showed that if $p_g$
is nonzero, then $A_0^0(F)$ is \lq\lq infinite-dimensional" [15], by showing 
that every component of $U^l\cap \pi_1^{-1}(X)$ has dimension at most $l$,
for any
$X$ in $S^lF$. From Roitman's theorems [17], one can deduce that $A_0^0(F)$
is finite-dimensional if and only if $\dim U^l_{se}\geq 3l$ for some $s,e,l$. 
So such a surface must have $p_g=0$, and the Bloch-Mumford conjecture
[4], states that
all surfaces with $p_g=0$ have finite-dimensional $A_0^0(F)$.
We refer to this as conjecture A. 

Roitman showed that the Albanese mapping  $\theta :A_0^0(F)\to
{\mathop{\rm Alb}}(F)$ is an isomorphism
on torsion. Since also the kernel of this map is divisible [6],
$A_0^0(F)$ is finite-dimensional, if and only if $\theta$ is an isomorphism.
 
Conjecture A  was shown to hold for surfaces with $p_g=0$ and  not of general 
type in [6], and for some examples of general type in (for example) 
[12, 3, 13, 25].
 
On the other hand, even surfaces with $p_g\not= 0$ have \lq\lq many" rational 
equivalences, for example coming from linear equivalences on curves on the 
surface, and moreover  
 Roitman's theorem \lq\lq closure of rational equivalence is 
albanese equivalence" [17]  implies that the 
Zariski closure of $U^l$ generates the Albanese kernel.
When $p_g \not= 0$ this implies that $U^l$ has infinitely many
components for large $l$.
Furthermore, a conjecture of Bloch and Beilinson [5] (which we refer to
as conjecture~B)  implies that if $F$ is
defined over $\bqb$, then $\theta$ is an isomorphism from $A_0^0(F_{\bqb})$ 
onto the $\bqb$ part of ${\mathop{\rm Alb}}(F)$. In this form, conjecture~B
is given as a question in~[17]. So far no surfaces with $p_g\neq0$
are known to satisfy conjecture B.

\mn{\sl4.2 Examples of  restatements  of conjectures\/}
 
\sn(1) General type surfaces over $\CC$.
Let $F$ be a surface of general type with $p_g=0$, with canonical 
divisor $K_F$. Then
conjecture A holds, \ie  $A_0^0(F) =0$, if and only if there exist 
integers $s,e$
such that for any pair of points $X$ and $Y$ on $F$ there exist effective
divisors $A,B$ in $ \vert sK_F\vert$ and $D$ in $\vert eK_F\vert$ with 
$X-Y= (A\cdot D)_F-(B\cdot D)_F +W$, where $W$ is a zero-cycle supported on 
the 
$-2$ curves of $F$ (and  $-1$ curves, if $F$ is non-minimal). Furthermore,
if $F$ belongs to a family, then $s$ and $e$ can be chosen to work
for all $F$ in the family.

\ni Proof. This follows from  {\sl3.10\/}\ (2), applied to the 
5-canonical model of $F$, say (in which case $s$ and $e$ are divisible by 5).
The map from $F$ to this model  contracts any $-1$ and $-2$ curves, but
is one-to-one elsewhere, when $F$ is of general type (see [7]). 

\sn(2) Surfaces over a number field. 
Let $F$ be a nonsingular surface in $\pt$ defined over $\bqb$. Then
conjecture B implies that $A_0^0(F_{\bqb}) =0$ 
(since $q=\dim{\mathop{\rm Alb}}(F)=0$), and  by {\sl2.1\/} this
 holds if and only if given points $X$ and $Y$ on $F_{\bqb}$ there are surface
 sections $A,B, D$ of $F$ such that $(A\cdot D)_F-(B\cdot D)_F=X-Y$.

If $F$ is any projective surface over $\bqb$,
conjecture B holds if and only if there
is a curve $E$ on $F$ such that given $X$ and $Y$ on $F$, there are
hypersurface sections $A,B,D$ and a zero cycle $Z$ supported on $E$
with $(A\cdot D)-(B\cdot D) + Z = X-Y$.

\mn{\sl 4.3 Algebraic formulation.}

\sn When the ideal defining $F$ is known, one can use (2.4) to
obtain other restatements.
For example, if $F$ is a surface given in $\pt$ by the homogeneous
 equation $f=0$, then two cycles $X,Y$ in $S^lF$
are rationally equivalent if and only if for some $s,e$ the equation     
$$
    \Phi_X  (g) R(b,h,f,g)= \Phi_Y  (g)R(a,h,f,g)\leqno(*)$$ 
for all $g$ in
the dual $\pt$, has a solution
$a,b,h$ (forms of degree $s,s,e$ on $\pt$) which is nontrivial
 (\ie neither side of this equation vanishes).
Here $R$ denotes a resultant of $4$ polynomials in $\pt$ (see [24] 
or [9]). Since for fixed $b,h,f$ the resultant $R(b,h,f,g)$
vanishes if and only if $g=0$ meets the intersection $\{b=h=f=0\}$,
it is the Chow polynomial (see [22]) for this intersection.
So conjecture B holds for nonsingular $F$ in $\pt$ defined over $\bqb$
if and only if equation (*) has nontrivial solutions for every pair of
points on $F$.
 
If $F$  is singular, one does not expect every pair of nonsingular 
points
to become rationally equivalent on the resolution. For example take $F$ to be 
a quartic nonsingular away from 
two disjoint double lines, whose resolution gives an elliptic ruled 
surface. But existence of solutions to (*) seems unlikely to depend on 
$F$ being nonsingular.If  solutions exist,they must either be
 forced to give cycles $A.D$,$B.D$ meeting the singular locus and
with high enough multiplicities there that the pull-backs can differ, or
else degenerate to trivial solutions.
For the above example this would give the correct result
$A^0_0(F' \backslash E)=0$,
where $F'$ is the normalization of $F$ and $E$ the elliptic curve over
which $F'$ is ruled (the double cover of one of the double lines on $F$).
This  suggests the following  version of conjecture B.

\sn{\sl Conjecture B. Let $F$ be any surface in $\pt$ defined over
$\bqb$, and $X,Y$ any two points in its nonsinular locus. Then equation (*) 
has nontrivial solutions.\/}

\sn This  with  {\sl3.8\/} and also {\sl4.4, 4.6\/} below suggest the
following questions.
{\sl
\sn(1) Let $F$ be any surface  in $\pt$ and  $X,Y$ any two points
on it. Is there  a sequence of triples of forms
$a_n,b_n,h_n$  of increasing degrees, constructed from the defining equations 
for $X,Y$, the equation for $F$, and their Taylor series,
such that whenever either $X$,$Y$ and $F$ are defined over $\bqb$ 
or $p_g(F)=0$, the triple
gives a solution to (*) for large $n$ (the value of $n$ depending
on the coordinates in the former case, see 4.4)?

\sn(2) Must there be such a sequence if the conjectures of [4] and [5] are 
true?\/}

\mn{\sl 4.4 Problem of finding  rational equivalences\/}

\sn We would like  to  see how to write down some explicit rational
equivalences on a surface $F$, and find $s,e$ as in {\sl3.10\/}\ (2)
(or {\sl4.2\/}\ (1) for general type) if $F$ is known to have $A_0^0(F)=0$.
For a trivial example, if $F$ is a plane and $X,Y$ two 
points on it, then we can write $X-Y=\bigl((L_1-L_2)\cdot L\bigr)$, 
where  $L$ is the
line joining $X$ to $Y$ and $L_1$, $L_2$ are other lines through $X$ and $Y$
respectively. So here $s=e=1$ will do. In {\sl4.6\/} we consider 
rational equivalences of points on surfaces in $\pt$, obtained by intersecting
with lines. Next one could look at the surfaces
treated in [6], elliptic surfaces with $p_g=0$. Here the Abel-Jacobi
theorem for curves (see [10]) is invoked, and the problem reduces to
writing down
explicit rational equivalences on a curve.   
For an example of general type, the generic Godeaux surface
$F=Q/{\ZZ_5}$, where $Q$ is a nonsingular quintic in $\pt$ on which
the cyclic group ${\ZZ_5}$ acts freely, was shown to have trivial Chow
group by Voisin [25]. Rational equivalences of points on $F$  
correspond to ${\ZZ_5}$-equivariant 
rational equivalences of orbits of points on the quintic, so this 
would also serve as an example  of the problem of  writing down some 
rational equivalences of cycles of low degree  on a general surface in $\pt$. 
This suggests the following.

\ni{\sl Conjecture: For a surface  $F$ in $\pn$, given $X,Y$ in
 $S^lF\times S^lF$ , 
there exist
 integers $s,e$ (given in terms  of $l$, invariants of $F$ and
its embedding, and of the coordinate field of $X$ and $Y$), such that 
$X\sim  Y$ if and only if $(X,Y)\in U_{e,s}^l$.}

 The existence of such $s,e$ is equivalent to the existence of
an upper  bound for the $s,e$ needed to express a rational equivalence,
by  {\sl3.4\/}. An analogous bound for divisors on 
curves, which is independent of coordinate fields, does exist.
For a nonsingular curve $C$ in $\pn$ of genus $g$ and degree $d$,
divisors $X$ and $Y$ (effective of degree $l$) are rationally equivalent
if and only if $X-Y=A-B$ for some sections $A,B$ of $C$ by hypersurfaces
of degree $s$, where Riemann-Roch gives that  any $s$ more than both
$(2g-2+l)/d$ and $r$ (where $r$ is the smallest for which the $r$-tuple
embedding of $C$ is projectively normal), will work.

  One consequence of the above conjecture is that  solving an equation
 similar to that in {\sl4.3\/}\  (given by {\sl2.4\/})  gives at least in 
principle a procedure for deciding 
whether or not $X\sim Y$ over a chosen algebraically closed field (or
whether $X-Y$ is torsion over a number field; see [4]),
 and for writing down a rational equivalence. Unfortunately it seems likely
that  the $s,e$  (if they exist) will be too large  for this to
be practical.
  
Also,  if conjecture B is assumed, then the existence of  bounds
dependent on coordinate field is equivalent to there being only
finitely many points on a surface over a number field (away from 
a finite collection of rational and elliptic
curves), when the surface has $p_g\neq 0$. This follows from Mumford's 
theorem
together with  Faltings' Mordell theorem 
(curves of genus at least $2$ have finitely
many rational points).

\mn{\sl4.5 Geometry  of the sets $U_{s,e}^l$.\/}

\sn  Two overlapping  approaches to these sets are to look directly at
the construction in {\sl3.2\/}, and to look at
 the determinantal equation  given by {\sl2.4\/} (as in {\sl4.2\/}).
This section consists of preliminary remarks on the former.


\ni(1) Expected role of $p_g$.

\ni To prove conjecture A in the form in {\sl4.1\/}, it would 
suffice to prove
that for a surface with $p_g=0$, there exist $s,e,l$ with
$\dim U_{s,e}^l\geq 3l$. This leads to the problem of finding $\dim\Delta'$
(where $\Delta'$ is the complement of the cover of the diagonal
in $(S^lF)^2$),  and
its fibre dimension over $U_{s,e}^l$. Equivalently, we can consider the image 
$T$ of $\Delta'$ in $(S^lF)^2\times S^{d-l}F$, and its projections
onto  $U_{s,e}^l$ and onto its image $S$ in  $S^{d-l}F$. Let
$\lambda$ be the fibre dimension of $T$ over $U_{s,e}^l$.  
Combining Bloch's conjecture and Mumford's theorem gives:

\ni{\sl Conjecture A: ${\dim}(T)-\lambda \geq 3l$ for some $s,e,l$ if and only
if $p_g=0$\/}.

\sn(2) Approximate calculations.

\ni The dimension of $T$ is at most $m+2k-x$, where 
$m=h^0(\co_F(e))-1$, $k \leq d-g$ and $g$ is the genus of generic section
of $F$ by an equation of degree $e$, and $x$ is the number of 
independent conditions imposed on a pair of complete intersections of 
type $s,e$ on $F$ with their degree $e$ equation in common,
by making them share $d-l$ points (counted with
 multiplicities).
For small $l$,  $T$ will for a general surface  be empty unless $S^{d-l}F$
 contains cycles
imposing fewer than $d-l$ conditions, and it is
hard to find $x$ when $T$ is nonempty without more knowledge of
$S$. On the other hand  for 
$l$ large (\ie greater than $g$), $U_{s,e}^l$ would be expected
to have  dimension less than
$3l$ (with exceptions for special surfaces such as ${\PP}^2$). 
 This suggests that we should first
try to analyse the set $S$  for small $l$ (see (3)). 
For a very rough check,
if we just assume $T$ is nonempty, then its dimension is at most 
$m+2k$, which by Riemann-Roch is linear in $\chi({\co_F})$. On the other hand,
using the nesting rules of {\sl2.5\/} and {\sl3.3\/},
 we can construct
whole families of rational equivalences $X\sim  Y$ from a given one, by 
adding ``redundant" expressions (so increasing $s$ and $e$). For this new
$s,e$, the fibre of $T$ over $X,Y$ has dimension depending on
at least  $\chi({\co_F})^2$. This  allows  the
 possibility of a role for
 $p_g(F)$  compatible 
 with  the above conjecture --- although these
 contributions might turn out to be insignificant.
 For a proof  we would need  to find more precise upper and lower bounds  
 for the fibre dimension, as well as bounds for $\dim T$. 

\sn(3) The set $S$

\ni For a given $l$, it is probably necessary to consider the sets 
 $S$ for  the case $l\leq g$ where $g$ is as in (2), and  because of 
``nesting''
({\sl3.3\/}) it is also sufficient. The set $S$ is then automatically the 
``special position locus'' for the image $J$
of $W_{s,e}$ in $S^{d-l}F$, which we define to be the set of cycles
$Z$ in $J$ imposing fewer conditions on complete intersections of type
$s,e$ than does the generic cycle in $J$.
In this range ($d-l$ large), it seems hard to apply intersection theory
(as in [26] or [10]) to find $S$ because it is an excess intersection 
and in general
will lie in the singular locus of $S^{d-l}F$.  

In [26], the collection of sets of $m$ points in special position
(in the ordinary sense of not spanning a ${\PP}^{m-1}$) 
on a general surface in ${\PP}^{3m-2}$ was shown to be finite, as
conjectured by Donaldson [8] in connection with Yang-Mills.
It would be interesting to know if
any similar conjecture could be made for the geometry of $S$.
A possible means of finding $S$ would be to try to generalize Serre's 
construction of rank 2 vector bundles associated to  special
position sets (where the scheme structure is also important; see [11]),
to obtain some moduli space of vector bundles
associated  to $S$ from which to estimate its dimension.

Often  $\dim S$ can be used to find $\dim T$ (so this would be
an alternative to finding the  $x$ in (2) above). 
For example, by increasing $s,e$ if necessary according to the rules of
{\sl2.5\/}, we may assume in addition to  
(i) $l\leq g$
as above, the condition
(ii) generic $Z$ in $S$ lies on a unique section $D$ of $F$ by an 
equation of degree $e$.
Then  $T\to S$ has infinite fibre over generic $Z$
if and only if some component of the corresponding $D$ has a linear system
of degree at most $l$. 
Choosing $s,e$ so that both $\dim S$ and $\dim T$ can be calculated
might present a problem, since the most natural candidate for a vector
bundle construction would be $s=e$. 
Also it may be that increasing $s$ and $e$
artificially conceals information ({\sl2.6\/} (1)).   
This applies particularly when it comes to estimating $\lambda$
(see (1)), unless 
the following holds.

\ni{\sl Conjecture: The component of $T$ with maximum fibre dimension
over $X,Y$ is that for which  the 
corresponding expression  $\bigl((A-B)\cdot D\bigr)_V=X-Y$ has $D=D_0+D_1$
with $(A\cdot D_0)_V=(B\cdot D_0)_V$ with $D_0$ of maximum degree.\/}

\sn(4) The case  $l=1$, $F$ of general type.

\ni If $F$ is a surface of general type  with $A_0^0(F)=0$, \ie 
satisfying conjecture A,
 then for  sufficiently large $s,e$ we have $U_{s,e}^1=F^2$.
Here the cover $T$ of  $S$
 must be generically finite on any  component of $T$ covering $F^2$
(since otherwise $F$ would be covered by rational curves). So
with the above notation, conjecture A can be stated as follows.

\sn{\sl Conjecture A: ${\dim}(S)-\lambda=4$, where $\lambda$ is the fibre
 dimension of $T$ over $F^2$.}

\ni The problem breaks into finding the two dimensions (or suitable
bounds). It might be hoped that the geometry of $S$ would give
information about $\lambda$.
Here  typical $Z$ in $S$ must determine a 
singular $D$, and have at least a triple point at a double point of $D$
 (\ie multiplicity higher than the multiplicity of the singularity on 
$D$).    

\eject

\ni{\sl4.6 Generic surfaces in $\pt$}

\sn(1) Proposition. {\sl The generic surface $F$ of degree $d\geq6$ in $\pt$ 
has $U^1_{1,1}\backslash\Delta=\emptyset$.}

\ni Proof. The set $U^1_{1,1}$ constructed in {\sl3.2\/} for $F\subset\pt$
is the set of pairs $(q_1,q_2)\in F\times F$ with $q_1-q_2=(L_1\cdot F)
-(L_2\cdot F)$ for some pair of lines $L_1$, $L_2$ in $\pt$. If
$F_d$ is the set of surfaces of degree $d$ in $\pt$, so $F_d$ is isomorphic
to the projective space of dimension ${d+3\choose3}-1$, it is easy to show 
that the correspondence $X_r$ below has dimension $\dim F_d-2r+8$, which
gives the result. The set $X_r$, a subset of $\pt\times{\mathop{\rm Gr}}(2,4)^2
\times F_d$, is defined to be the closure of the set of all $(p,L_1,L_2,F)$
such that $L_i\not\subset F$, $L_1\neq L_2$, $(L_i\cdot F)\supseteq rp$.

\ni Using the Taylor series for $F$ at $p$ (see (3) below), we can prove the
following.

\sn(2) Proposition. {\sl Let $U$ be the closure of 
$U^1_{1,1}\backslash\Delta$. For generic quintic $F$, $U$ is finite (and 
nonempty). For generic quartic
$F$, $U$ is a surface with projection onto $F$ of degree 24. For $d\leq3$,
$U=F\times F$ (hence $A^0_0(F)=0$, as is well-known).}

\ni Corollary. {\sl For quartic $F$, where $U^1=\{(q_1,q_2)\in F\times F:
q_1\sim q_2\}$ as in 3.1, we have $\overline{U^1}=F\times F$.}

\ni To prove this, we consider the surface $U$ above and its 
\lq\lq iterates\rq\rq, \ie sets of pairs $(q_1,q_2)\in F\times F$
such that $q_1-q_2=\sum_{i=1}^k\bigl((L_{i1}-L_{i2})\cdot F\bigl)$ for
$k$ pairs of lines. This also verifies Roitman's theorem on closure of
rational equivalence in this case. Also, $U^1$ contains 
self-products of infinite sequences of curves, derived from \eg 
(a) rational
curves on $F$, and (b) the ``touch'' curve $R$ on $F$ (\ie the
curve of points $p$ where some line $L$ gives $(L\cdot F)=4p$). Starting
with a curve $C\subset F$ of points rationally equivalent to each other, and
adding $\pi_1\bigl(\pi_2^{-1}(C)\cap U^1_{1,1}\bigr)=C'$, gives
$(C'\cup C)\times (C'\cup C)\subset U^1$. This construction is then repeated. 
This supports conjecture B. It would be interesting to know if the union of 
$U$ and its iterates should  contain all of
$F_{\bqb}\times F_{\bqb}$, when $F$ is a quartic defined over a number
 field satisfying conjecture B.

\sn(3) Taylor series. Let $p$ be a point on the surface $F$ in $\pt$.
Let $f=\sum(f_i)_p$ be the Taylor expansion for a polynomial $f$ defining
$F$ at $p$. If coordinates $(X,Y,Z,T)$ are chosen for $\pt$, and
$p\in\{T\neq0\}$, then with respect to $x=X/T$, $y=Y/T$ and $z=Z/T$,
we have
$$
(f_i)_p(a,b,c)=\sum_{\scriptstyle l+m+n=i\atop\scriptstyle l,m,n\geq0}
 {\pa^i f \over \pa x^i\, \pa y^m\, \pa z^n}\bigg\vert_p
 {a^l b^m c^n\over l!\, m!\, n!}.
$$

\sn{\sl Lemma. For $p\in F$, we have $(f_0)_p=0$. The map $L\mapsto
L\cap\{T=0\}$ is a one-to-one correspondence between the set of lines
$L$ in $\pt$ such that $(L\cdot F)\supset rp$, and the set of points
$(a,b,c,0)$ in $\pt$ with $(a,b,c)$ a solution of $(f_1)_p=\ldots=
(f_{r-1})_p=0$.}

\sn Definition. The ``polar locus'' $F^r_q$ for a point $q$ (which may be on 
$F$) is the set of points $p$ in $F$ for which there is some line $L$
giving $(L\cdot F)\supseteq rp+q$.

\ni This generalizes the classical polar locus (which is $F^2_q$).

\sn{\sl Corollary of Lemma. If $q=(a,b,c,0)$, then
$F^r_q$ is the set of points $p\in F$ such that
$(f_1)_p(a,b,c)=\ldots=(f_{r-1})_p(a,b,c)=0$.}

\mn(4) Examples and proof of (2). If $q$ lies on a quartic $F$, then $F^3_q$
consists of 24 points in general. At a generic point $p_i\in F$, there are
two lines $L_1$ and $L_2$ for which $(L_i\cdot F)\supset3p_i$ (by the lemma);
so if $p_i$ is one of the $24$ points in $ F^3_q$, one of these lines 
(say $L_1$) has $(L_1\cdot F)=3p_i+q$.
Then $(L_2\cdot F)=3p_i+q_i$, some $q_i$. This gives rise to 24 points $q_i$
rationally equiavalent to $q$ (as in {\sl4.6\/}\ (2)).

To prove that $U$ in {\sl4.6\/}\ (2) is finite for a generic quintic, a
dimension count shows that it is enough to find some nonsingular quintic
for which $U$ is nonempty. The generic quintic $F$ defined on $\{T\neq0\}$
by $f\in{\cal L}$ below works:
$$
  {\cal L} = \langle x+y+z, xy, xyz, f_4, f_5 \rangle,
$$
with $f_4$ and $f_5$ generic forms of degree 4 and 5. For let $p=(0,0,0,1)$,
$r=(0,-1,1,0)$, $s=(1,0,-1,0)$, $L_1=\langle p,r\rangle$ and
$L_2=\langle s,r\rangle$. Then $(L_i\cdot F)\supset 4p$, but not $5p$ in
general (for $i=1,2$), so $U\neq\emptyset$. By Bertini's theorem, $F$ is
nonsingular, since the base locus of ${\cal L}$  on $\pt$ consists of $p$, 
and generic $F$ is nonsingular at~$p$.
\bigskip
\ni {\bf Acknowledgements.}
The author is grateful to
the Max Planck Institut in Bonn, where this work was begun;
and thanks the referee for helpful comments. In particular, the referee
pointed out that in the proof of {\sl 2.3\/} the underlying field need only
be  infinite and perfect, and suggested a different method for finite fields.

\bigbreak
\baselineskip=3mm

\ni{\bf References.}
\medskip
\item{1.} Akizuki, Y.: Theorems of Bertini on linear systems. J. Math.
Soc. Jpn. {\bf3}, 170--180 (1951)
\item{2.} Baldassari, M.: Algebraic Varieties. (Ergebnisse der Mathematik)
Berlin Heidelberg New York: Springer 1956
\item{3.} Barlow, R.: Rational equivalence of zero-cycles for some more 
surfaces with $p_g=0$. Invent. Math. {\bf79}, 303--308 (1985)
\item{4.} Bloch, S.: Lectures on Algebraic Cycles. (Duke University Math.
Series IV) Durham, NC: Duke University 1980
\item{5.} Bloch, S.: Algebraic cycles and values of L-functions. J.
Reine Angew. Math. {\bf350}, 94--108 (1984)
\item{6.} Bloch, S., Kas, A. and Lieberman, D.: Zero cycles on surfaces with 
$p_g=0$. Compositio Math. {\bf33}, 135--145 (1976)
\item{7.} Bombieri, E.: Canonical models of surfaces of general type.
Publ. Math IHES {\bf42}, 171--220 (1973)
\item{8.} Donaldson, S. K.: Instantons in Yang-Mills theory. In:
Interactions between particle physics and mathematics, pp. 59--75.
Oxford: Oxford University Press 1989
\item{9.} Fulton, W.: Intersection Theory. (Ergebnisse der Mathematik)
Berlin Heidelberg New York: Springer 1984
\item{10.} Griffiths, P.A.: An introduction to the theory of special divisors 
on  algebraic curves. (CBMS regional conference series 44)
Providence: American Mathematical Society 1980
\item{11.} Griffiths, P. and Harris, J.: Residues and zero-cycles on algebraic
varieties. Ann. Math. {\bf108}, 461--505 (1978)
\item{12.} Inose, H. and Mizukami, M.: Rational equivalence of zero-cycles 
on some surfaces with $p_g=0$. Math. Ann. {\bf244}, 205--217 (1979)
\item{13.} Keum, J.H.: Some new surfaces of general type with $p_g=0$. 
Preprint, University of Utah, Salt Lake City.
\item{14.} Lewis, J.D.: A Survey of the Hodge Conjecture. Montreal:
CRM publications 1991
\item{15.} Mumford, D.: Rational equivalence of 0-cycles on surfaces.
J. Math. Kyoto Univ. {\bf9}, 195--204 (1969)
\item{16.} Murthy, M. and Swan, R.: Vector bundles over affine surfaces.
Invent. Math. {\bf 36}, 125--165 (1976)
\item{17.} Roitman, A.A.: The torsion of the group of zero-cycles modulo 
rational equivalence. Ann. Math. {\bf111}, 553--569 (1980)
\item{18.} Samuel, P.: M\'ethodes d'Alg\`ebre Abstraite en G\'eom\'etrie 
Alg\'ebrique. Berlin Heidelberg New York: Springer 1955
\item{19.} Samuel, P.: Rational equivalence of arbitrary cycles. Am. J.
Math. {\bf78}, 383--400 (1956)
\item{20.} Severi, F.: Un'altra propriet\`a fondamentale delle serie di 
equivalenza sopra una superficia. 
Rend. Accad. Linc. {\bf21}, 3--7 (1935)
\item{21.} Severi, F.: Serie, sistemi d'equivalenza, e corrispondenze 
algebriche sulle variet\`a algebriche (vol 1). Rome: Cremonese 1942
\item{22.} Shafarevich, I. R.: Basic Algebraic Geometry.
Berlin Heidelberg New York: Springer 1974
\item{23.} Todd, J. A.: Some group-theoretic considerations in algebraic 
geometry, Ann. Math. {\bf35}, 702--704 (1934)
\item{24.} Van der Waerden, B.L.: Modern Algebra, volume II. New York:
Unger 1950
\item{25.} Voisin, C.: Sur les z\'ero-cycles de certaines 
hypersurfaces munies d'un automorphisme. Ann. Scuola Norm. Sup. Pisa Cl.
Sci. (4) {\bf19}, 473--492 (1992)
\item{26.} Xu, M.W.: The configuration of a finite set on surface.
Preprint, MPI Bonn (1990)
\item{27.} Zariski, O.: Algebraic Surfaces.
(Ergebnisse der Mathematik) Berlin Heidelberg New York:
Springer 1971

\bye